\documentclass[a4paper,english]{article}

\usepackage[english]{babel}
\usepackage[utf8]{inputenc}
\usepackage{hyperref}
\usepackage{amsmath}
\usepackage{amsfonts}
\usepackage{amssymb}
\usepackage{mathtools}
\usepackage{amsthm}
\usepackage{graphicx}
\usepackage{tikz}
\usepackage{times}
\usepackage{graphicx}
\usepackage{array}
\usepackage{amsmath}
\usepackage{amsfonts}
\usepackage{amssymb}
\usepackage{mathtools}
\usepackage{enumitem}
\usepackage{todo}
\usepackage{gensymb}
\usepackage{xcolor}
\usepackage{booktabs}
\usepackage{colortbl}
\usepackage{makecell}
\usepackage{hyperref}
\usepackage[dvipsnames]{xcolor}
\definecolor{ecosblue}{RGB}{49, 123, 188}

\usetikzlibrary{patterns, shapes, arrows, positioning, fit, backgrounds, calc}

\usepackage{zibtitlepage}
\usepackage{authblk}

\definecolor{myblue}{RGB}{7, 127, 247}
\definecolor{myred}{RGB}{203, 65, 84}
\definecolor{mygreen}{RGB}{85, 107, 47}

\usepackage[a4paper, total={6in, 9in}]{geometry}

\begin{document}

\title{Benchmarking Realistic Synthetic Instances Against a Large-Scale District Heating Network: A Multi-Objective Optimization Study for Berlin}

\author[1]{Annika Buchholz\footnote{corresponding author, buchholz@zib.de, \ZTPOrcid{0009-0008-3235-0896}}}
\author[1]{Stephanie Riedmüller\footnote{riedmueller@zib.de, \ZTPOrcid{0009-0006-4508-4262}}}
\author[1]{Matthew Passage\footnote{passage@zib.de}}
\author[1]{Janina Zittel\footnote{zittel@zib.de, \ZTPOrcid{0000-0002-0731-0314}}}

\affil[1]{Zuse Institute Berlin, Applied Optimization, Berlin, Germany}

\maketitle

\begin{abstract}

Decarbonizing urban energy systems requires optimization approaches capable of handling the operational complexity of large-scale district heating networks. However, existing studies typically focus on a single, real-world network, which makes results difficult to compare and limits the transferability of insights to other systems. To overcome these challenges and enable a more systematic evaluation, realistic synthetic instances play an important role: they provide controlled, reproducible environments for testing optimization algorithms independent of a specific case study while still capturing essential structural and temporal characteristics of real systems. When designed carefully, such instances enable systematic benchmarking, facilitate methodological development, and support comparative studies across different algorithms and modeling choices.

In this work, we generate a suite of large-scale synthetic instances for multi-objective optimization of district heating systems. The instances are openly available in the form of the underlying network topology in a JSON format and also as a mixed integer program (MIP) as MPS files to enable modeling and algorithmic benchmarking. They are constructed using a transparent procedure that allows to recreate or extend the given instances or apply the same procedure to other network based problems. We apply this methodology to Berlin’s district heating network, the most complex in Western Europe, formulating a tri-objective mixed-integer model focused on unit commitment over a horizon of up to 25 years with a 4-hour temporal resolution. The corresponding real-world instance serves as a reference point to demonstrate that the generated datasets reflect key behaviors of actual district heating operations.

A computational study provides a detailed comparison between the synthetic instances and the real-world Berlin data, showing under which conditions the generated instances reproduce realistic optimization characteristics. Furthermore, we investigate which features make the resulting models computationally challenging. The findings highlight how well-designed synthetic instances can support robust benchmarking practices and enable meaningful assessment of (multi-objective) optimization methods for large-scale district heating systems.
\end{abstract}

\section{Introduction}
While district heating is recognized as a cornerstone of the thermal transition, the field lacks a standardized, scalable benchmark for evaluating optimization methods across diverse network topologies and investment horizons. The highly complex mixed-integer optimization problems involved in planning and designing large-scale district heating networks feature strong temporal coupling and heterogeneous system components. Such models quickly become computationally challenging when integrating operational, structural, and investment decisions \cite{jerry}. Despite methodological advances, empirical evaluation remains limited: most studies rely on single case studies based on proprietary data and simplified settings, restricting reproducibility and cross-study comparability \cite{Sporleder2022}. At the same time, increasing model detail through sector coupling, high-resolution temporal modeling, and multi-objective formulations further exacerbates computational difficulty \cite{guo2024novel}. Due to a lack of established benchmark libraries, systematic benchmarking for district heating and multi-energy systems is limited, hindering the analysis of structural drivers of problem difficulty \cite{PereaMoreno2026, jerry}. To address this gap, we generate structured synthetic instance families derived from a real-world reference system in Berlin. By varying structural and parametric properties while preserving key system characteristics, the approach enables controlled, reproducible, and realistic computational experiments. Our main contributions are 
\begin{itemize}
    \item the methodology to generate a meaningful set of instances based on a given instance of a multi-energy network, that vary a predefined set of parameters and structural properties, while maintaining key features,
    \item a computational study based on the district heating network of Berlin, including an evaluation of which properties make the problem difficult,
    \item the publication of the network topology of the generated instances and mixed-integer programs for a corresponding unit commitment problem (Zenodo project 19550701 \cite{Zenodo}), together with the open-source generator code (GitLab project multi-energy-instance-generator \cite{generator}).
\end{itemize} 

To support diverse research objectives, we provide network topologies that remain agnostic to specific model constraints; meanwhile, the MPS files represent a finalized optimization instance of the DHOT model \cite{riedmuller2025enhancing}, serving as a standardized testbed for benchmarking advanced solution heuristics.


\section{Model base} \label{Introduction} 
\begin{figure}[h!]
    \centering
    \resizebox{1.0\linewidth}{!}{
        \usetikzlibrary{positioning,calc,fit} 
\tikzstyle{container} = [rectangle, minimum width=1.5cm, minimum height=0.6cm, text centered, draw=BrickRed, fill=BrickRed!30] 

\tikzstyle{containergreen} = [rectangle, minimum width=1.5cm, minimum height=0.6cm, text width=2.2cm, text centered, draw=OliveGreen, fill=OliveGreen!30] 

\tikzstyle{container1} = [rectangle, minimum width=1.5cm, minimum height=0.6cm, text centered, draw=lightgray, fill=lightgray] 

\tikzstyle{container2} = [rectangle, minimum width=1.5cm, minimum height=0.6cm, text width=2.2cm, text centered, draw=ecosblue, fill=ecosblue!30] 

\tikzstyle{container3} = [rectangle, minimum width=1.5cm, minimum height=0.6cm, text centered, draw=BrickRed] 

\tikzstyle{containeremission} = [rectangle, minimum width=1.5cm, minimum height=0.6cm, text centered, draw=black] 

\tikzstyle{oval} = [rectangle, minimum width=1.5cm, minimum height=0.8cm, text centered, dashed, draw=BrickRed, rounded corners=5pt] 

\tikzstyle{ovalgreen} = [rectangle, minimum width=1.5cm, minimum height=0.8cm, text centered, dashed, draw=OliveGreen, rounded corners=5pt] 

\tikzstyle{site} = [draw=black, dashed, rounded corners, inner sep=0.5cm] 

\tikzstyle{arrow} = [thick,->,>=stealth] 

\begin{tikzpicture}[node distance=1.2 cm] 
\node (demand1) [container] {Demand 1}; 
\node (hbal1) [oval, below =4.5cm of demand1] {Heat Balance}; 
\node (pbal1) [ovalgreen, right=2cm of hbal1] {Power Balance}; 
\node (H1) [container1, below of=hbal1] {Heating Plant}; 
\node (CHP1) [container1, right=1.5cm of H1] {CHP 1}; 
\node (demand2) [container, right=7cm of demand1] {Demand 2}; \node (hbal2) [oval, below =4.5cm of demand2] {Heat Balance}; \node (pbal2) [ovalgreen, right=2cm of hbal2] {Power Balance}; \node (CHP2) [container1, below of=hbal2] {CHP 2}; \node (CHP3) [container1, right=1.5cm of CHP2] {CHP 3}; 
\node (pump1) [container3, above=1.2cm of hbal1] {Pump}; \node (pump2) [container3, above=1.2cm of hbal2] {Pump}; 
\node (hbal3) [oval, below =1cm of demand1] {Heat Balance}; \node (hbal4) [oval, below =1cm of demand2] {Heat Balance}; \node (pump3) [container3, right=2.5cm of hbal3] {Pump}; 
\node (storage) [container3, left=1.0cm of hbal3] {Storage}; 
\node (p2h) [container1, right=1.5cm of CHP3] {Power to Heat}; 
\node (fuel1) [container2, below=1cm of H1] {Fuel Market \\ Gas}; 
\node (fuel2) [container2, right=1.5cm of fuel1] {Fuel Market \\ Syn Gas}; 
\node (fuel3) [container2, right=1.5cm of fuel2] {Fuel Market \\ Biomethane}; 
\node (fuel4) [container2, right=1.5cm of fuel3] {Fuel Market \\ Biomass}; 
\node (emission) [containeremission, left=1cm of fuel1] {Emission}; 
\node (pimport) [containergreen, below=1cm of p2h] {Power Import}; 
\node (power) [containergreen, above=3.3cm of p2h] {Power Market}; 
\node (site1) [site, fit=(H1)(CHP1)(hbal1)(pbal1), label={[above right, xshift=18mm]Prod. Site 1}] {}; 
\node (site2) [site, fit=(CHP2)(CHP3)(hbal2)(pbal2)(p2h), label={[above right, xshift=28mm]Prod. Site 2}] {}; 

\draw [arrow, BrickRed] (pump1) -- (hbal3); \draw [arrow, BrickRed] (hbal3) -- (demand1); \draw [arrow, BrickRed] (H1) -- (hbal1); 
\draw [arrow, BrickRed] (CHP1) to[out=140, in=0] (hbal1); 
\draw [arrow, BrickRed] (pump2) -- (hbal4); \draw [arrow, BrickRed] (hbal4) -- (demand2); \draw [arrow, BrickRed] (CHP2) -- (hbal2); \draw [arrow, BrickRed] (CHP3) -- (hbal2); \draw [<->, BrickRed] (hbal3) -- (pump3); 
\draw [<->, BrickRed] (pump3) -- (hbal4); 

\draw [arrow, BrickRed] (hbal1) -- (pump1); 
\draw [arrow, BrickRed] (hbal2) -- (pump2); 
\draw [arrow, OliveGreen] (CHP1) -- (pbal1); 
\draw [arrow, OliveGreen] (CHP2) -- (pbal2); 
\draw [arrow, OliveGreen] (CHP3) -- (pbal2); 
\draw [<->, OliveGreen] (pbal2) to[out=70, in=260] (power); 
\draw [<->, OliveGreen] (pbal1) to[out=30, in=-150] (power); 
\draw [arrow, OliveGreen] (pbal1) to[out=150, in=350] (pump1); 
\draw [arrow, OliveGreen] (pbal2) to[out=150, in=350] (pump2); 
\draw [arrow, OliveGreen] (pimport) -- (p2h); 
\draw [arrow, OliveGreen] (pbal1) to[out=120, in=250] (pump3); 
\draw [arrow, OliveGreen] (pbal2) to[out=120, in=-20] (pump3); 
\draw [arrow, BrickRed] (p2h) to[out=180, in=0] (hbal2); 
\draw [arrow, ecosblue] (fuel1) -- (H1); 
\draw [arrow, ecosblue] (fuel1) to[out=30, in=200] (CHP1); 
\draw [arrow, ecosblue] (fuel2) to[out=90, in=270] (CHP1); 
\draw [arrow, ecosblue] (fuel3) to[out=70, in=270] (CHP2); 
\draw [arrow, ecosblue] (fuel4) -- (CHP3); 
\draw [arrow, ecosblue] (fuel1) to[out=20, in=190] (CHP2); 
\draw [arrow, black] (fuel1) -- (emission); 
\draw [<->, BrickRed] (hbal3) -- (storage); 
\end{tikzpicture}
    }
    \caption{Simplified illustration of an instance with two production sites (fuel (blue), heat (red), power (green)), emission (black).}
    \label{fig:MV}  
\end{figure}
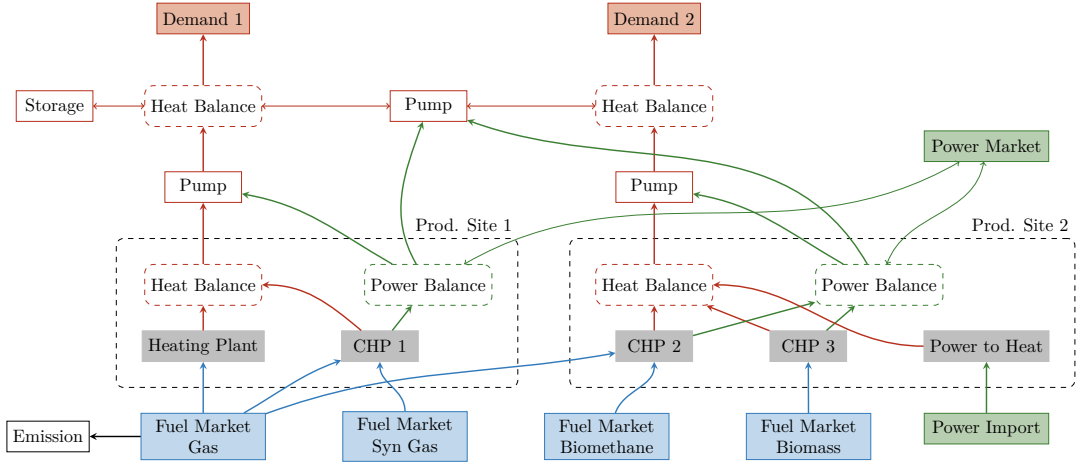

In the following, we consider multi-energy systems based on a network structure. The system is represented as a directed, unweighted graph $G=(V,A)$, where $V$ denotes the set of nodes and $A$ the set of directed arcs connecting these nodes. The arcs represent feasible flows of resources such as heat, fuel, power, steam, and emissions between nodes, for example, via pipelines or transmission lines. 
The nodes correspond to system components (assets) that fulfill different functional roles: they may convert one set of resources into another (converters), store a resource (storage), enable the import or export of a resource (markets), collect and distribute flows (balance nodes), or specify the demand for a given resource.
Fig.~\ref{fig:MV} illustrates such a network structure using the example of two production sites of the district heating network.

\section{Instance generation}
We aim to design the framework in a way that enables systematic scalability, controlled structural variation, and realistic parametrization. To achieve this, we introduce a transparent and reproducible procedure for generating realistic synthetic multi-energy instances. This process defines component configurations, constructs network topologies, generates time-series data, and exports the resulting optimization models in standardized formats.

\subsection{Generator parameters}

The generator is controlled by configurable global and temporal parameters, network component parameters, structural parameters, and stochastic and scenario parameters, summarized in Table \ref{tab:parameters}.


\begin{table}[h]
    \centering
        \caption{Parameters of instance generation}

    \begin{tabular}{
    >{\raggedright\arraybackslash}p{3cm} 
    >{\raggedright\arraybackslash}p{3cm} 
    >{\raggedright\arraybackslash}p{8 cm}
}
         \hline
         \textbf{Category} &	\textbf{Parameter Group} &	\textbf{Key Parameters \& Features} \\
         \hline
         Global \& Temporal	& Scope \& Economics	& Planning horizon, temporal resolution, constant discount rate, and constant inflation rate.\\
         \hline
         Network components	& Demand &	Number of nodes, associated heat demand time series, transport capacity/costs between heat demand nodes. \\
         & Storage	& Number of thermal storage units, connected to a local heat balance node through inflow/outflow heat arcs. \\
         & Converters	& Number of converters, fuel-to-heat, combined heat and power plant (CHP), and power-to-heat units, heat pumps characterized either by a fixed conversion efficiency or by a piecewise-linear conversion function, each requires one or two specific input fuels. \\
        & Markets	& Number of fuel supply options (gas, oil, coal, biomass, etc.) with distinct price series, emission certificates, electricity import/export nodes. \\
\hline
     Structural & Production Sites & Number of sites (urban districts), each conversion unit is assigned to exactly one site with local heat/power balance nodes. \\
     & Fuel \& Power Transport & Stochastic capacity constraints ($\lambda_{fuel}$), stochastic transport costs ($\kappa_{fuel}$). \\
     & Heat Transport & Site-to-demand connectivity, heat transport nodes, stochastic capacity ($\lambda_{heat}$), and power-driven transport ($\kappa_{heat}$). \\
     \hline
     Stochastic \& Scenario & Variability & Conversion efficiencies, storage losses, transport costs, capacity levels, and ramping limits.\\
     & Time Series & Emission factors/prices, synthetic temperature, demand, fuel, and power price profiles. \\
     \hline
      
    \end{tabular}
    \label{tab:parameters}
\end{table}

\subsection{Generation of components} \label{comp gen}

All technical parameters for conversion units and storage units, including capacities, efficiencies, losses, ramping limits, minimum up- and down-times, and start-up costs, are generated independently for each node during instance generation by sampling from uniform distributions over predefined parameter intervals. These parameters remain constant throughout the planning horizon.

The conversion of fuel to heat, fuel to heat and power, and power to heat is represented by characteristic curves defined by a source series $y_1$ and one or two target series $y_2$. These series specify the defining points of the piecewise-linear input-output relation. For one-to-one converters (e.g., fuel to heat and power to heat), the curve is constructed by using the beforehand randomly sampled minimum and maximum output capacities $q_{\min}$ and $q_{\max}$. These capacities are scaled by random factors drawn from a uniform distribution $\mathcal{U}$ of a predefined interval $\alpha_{\min}, \alpha_{\max} \sim \mathcal{U}(\underline{\alpha},\overline{\alpha})$, resulting in
\[
y_1 = \left[\frac{q_{\min}}{\alpha_{\min}}, \frac{q_{\max}}{\alpha_{\max}}\right], 
\qquad
y_2 = \left[q_{\min}, q_{\max}\right].
\]

For one-to-two converters (e.g., combined heat and power units), the second target series is generated by scaling $q_{\min}$ and $q_{\max}$ with two factors drawn from distinct distributions $\beta_{\min} \sim \mathcal{U}(\underline{\beta_{\min}},\overline{\beta_{\min}})$ and $\beta_{\max} \sim \mathcal{U}(\underline{\beta_{\max}},\overline{\beta_{\max}})$. For a user-defined subset of converters, conversion ratios are used instead of characteristic curves. The conversion ratio $\rho$ is sampled independently from a predefined interval chosen to respect physical laws, $\rho \sim \mathcal{U}(\underline{\rho},\,\overline{\rho})$. The conversion is modeled as $x^{\text{out}}_{t} = \rho \, x^{\text{in}}_{t}$, where $x_t$ denotes the flow variable at time $t \in T$.
To ensure feasibility, the total heat production capacity is enforced to be sufficient to cover the aggregated heat demand.

Emissions are modeled using an emission factor $r_{f}$ sampled from $r_{f} \sim \mathcal{U}(\underline{r_f},\overline{r_f})$. Emissions are linked to an emission certificate node with a certificate price $r_c \sim \mathcal{U}(\underline{r_c},\overline{r_c})$, resulting in emission costs $r_f \cdot x_{\text{fuel}, t} \cdot r_c$, where $x_{\text{fuel},t}$ denotes the consumed fuel at time $t\in T$. Transport costs are modeled as nodes with technology-specific characteristic curves that map transported fuel flow (source series) to cost (target series), generated from predefined intervals. Transportation capacity limits are likewise represented as nodes, using values drawn from a uniform distribution over a specified interval. For the transport of heat, electric pumps are required. They follow a similar approach to converters, with characteristic curves defined analogously but using distinct parameter ranges.

\subsection{Construction of network topology} \label{constrction}

The topology is generated algorithmically using the following procedure.

\begin{enumerate}
    
\item 
Converters are assigned to production groups by uniform random sampling. Consequently, the number of converters per group may vary substantially, reflecting the heterogeneous structure of real-world energy systems.

\item 
Each production group is connected randomly to a demand node via heat balance nodes. 

\item The heat balance nodes associated with a given demand are interconnected through transportation cost nodes and pumps. With probability $\lambda_{heat} \in [0,1]$, a capacity constraint on a heat link is added, with probability $\kappa_{heat} \in [0,1]$ the heat transport requires electricity input.

\item 
Each storage unit is connected to exactly one heat balance node through two directed edges representing loading and unloading.

\item 
Converters are connected to at most two fuel market nodes via transportation or 
capacity-limiting utility nodes, reflecting fuel specificity. Market nodes $M$ are 
ordered by fuel type (natural gas, synthetic gas, oil, coal, biomethane, biomass) and 
assigned strictly decreasing weights $w_i = |M| - i$, $i = 0, \dots, |M| - 1$. These weights bias the assignment of market–converter edges such that higher-weight fuels (like natural gas) are more likely to be selected as inputs, while lower-weight fuels (like biomass) are chosen less frequently. For each converter $c \in C$, a pair of market nodes $m_i, m_j \in M$, is selected 
greedily: at each step the market node most under-represented relative to its 
weight-proportional edge budget is preferred, ranking all nodes by the deficit $\delta_i^+ -  \frac{2|C|\, w_i}{\sum_j w_j}$
where $\delta_i^+$ denotes the number of edges already assigned to $m_i$. The two 
lowest-ranked nodes are assigned to $c$.



\item With probability $\lambda_{fuel} \in [0,1]$, a capacity constraint modeled by an edge utility node is added between a fuel import market and a converter. With probability $\kappa_{fuel} \in [0,1]$, a transport cost component is added.

\item Each fuel that produces emissions is connected to the CO$_2$ market node.

\item Each converter that transforms electricity into heat is connected to a local power balance node, which in turn is linked to the power market via incoming and outgoing edges, while each pump node receives an incoming edge from a power balance node.

\end{enumerate}

Direct coupling of converters (e.g., fuel-to-steam followed by steam-to-heat) is not modeled explicitly. Instead, such chains are represented by a single conversion mapping. For example, a fuel-to-steam conversion $\varphi$ followed by steam-to-heat conversion $\vartheta$ is contracted into a fuel-to-heat conversion $\vartheta \circ \varphi$ by a composition of their corresponding characteristic curves and conversion ratios.

\subsection{Generation of temporal profiles} \label{gen temp}
Temporal profiles are generated using methods tailored to their respective applications, as different types of time series exhibit distinct statistical properties and dynamics.

\subsubsection{Temperature forecast profiles}

Temperature profiles are generated by combining a seasonal baseline with stylized diurnal fluctuations, a technique that aligns with existing energy demand and meteorological models employing periodic harmonic components and stochastic noise \cite{wilks1999weather, horsch2018pypsa}. Let each time step $t\in T$ be associated with a month $s(t)\in S = \{1, \dots , 12 \}$. For each month, a fixed average temperature level $\bar{\theta}_{s(t)}$ is specified. The temperature at time step $t$ is then constructed as
\[
\theta_t
=
\bar{\theta}_{s(t)}
+
\delta_t,
\]
where the deviation term $\delta_t$ represents stylized intra-day fluctuations. Given the four-hour temporal resolution, one day is represented by six consecutive time steps. Thus,
\(
k(t) = (t-1) \bmod 6
\)
denotes the position of time step $t$ within the corresponding day. The intra-day fluctuation is modeled as
\[
\delta_t
=
\xi_t \, \bigl(-\cos(\tfrac{\pi}{3} k(t))\bigr),
\]
where the amplitude factor $\xi_t$ is drawn independently for each time step from a uniform distribution,
\(
\xi_t \sim \mathcal{U}(-4,4).
\)
This construction produces a repeating sinusoidal temperature pattern over each day, centered around the monthly mean temperature and perturbed by random amplitude variations. The resulting time series, therefore, captures both seasonal temperature differences and stylized diurnal variability, while remaining lightweight and fully reproducible.

\subsubsection{Heat demand profiles}
Heat demand profiles are generated in two stages. First, a single aggregate heat demand profile $d$ for the entire system is generated using a seasonally structured, correlated stochastic process. In a second step, this aggregate profile is split into multiple spatially distributed demand profiles $d_{i,t}$ for each demand node $i \in D$ representing the urban districts and each time step $t \in T$. This autoregressive structure captures short-term persistence and seasonal behavior, which is consistent with existing heat load modeling and forecasting literature \cite{heatdemand}.

For each month $s$, we specify a demand interval
$[\underline d_{s},\,\overline d_{s}]$, using the real-world Berlin model as reference. From this interval, a segment-specific reference level $\mu_{s} = \frac{\underline d_{s} + \overline d_{s}}{2}$ and dispersion $\sigma_{s} = \frac{\overline d_{s} - \underline d_{s}}{6}$ are derived.
The raw aggregate demand trajectory is generated by a seasonally varying first-order autoregressive process. For all time steps $t \geq 2$,
\[
 d_{t}
=
\mu_{s(t)}
+
\phi\bigl( d_{t-1}-\mu_{s(t-1)}\bigr)
+
\varepsilon_t,
\qquad
\varepsilon_t \sim \mathcal{N}\!\left(0,\sigma_{s(t)}^2\right),
\]
where $\phi\in(0,1)$ is a fixed autocorrelation parameter controlling the smoothness of the profile and $\varepsilon_t$ denotes the noise. The initial value $ d_{1}$ is drawn uniformly from the admissible interval of the first segment. To ensure physically meaningful values and to strictly respect the user-defined seasonal ranges, all generated values are clipped to the admissible interval, $d_{t} =
\min\{\overline d_{s},\,
\max\{\underline d_{s}, d_{t}\}\}$, for $t$ in month $s$. 

This aggregate profile $d_t$ is subsequently split into $|D|$ individual demand profiles. To this end, a vector of structural shares $(\omega_1,\ldots,\omega_{|D|})$ is sampled from a symmetric Dirichlet distribution with all parameters fixed to $0.3$, yielding persistent heterogeneous demand sizes across nodes.

Small time-dependent perturbations are introduced to avoid perfectly fixed proportional allocations. For each time step $t$ and each demand node $i$, a multiplicative noise factor $\eta_{i,t}$ is drawn independently from a Gamma distribution with parameter $50.0$ and scale $1$. The time-dependent shares are then computed as
\[
\omega_{i,t}
=
\frac{\omega_i\,\eta_{i,t}}{\sum_{j=1}^{|D|} \omega_j\,\eta_{j,t}}.
\]

Finally, the individual demand profiles are obtained by
\[
d_{i,t} = \omega_{i,t}\, d_t,
\qquad
i \in D,\;\; t\in T.
\]

This procedure yields structurally persistent large and small demand nodes, representing different districts, while still allowing for moderate temporal variability in their relative share.

\subsubsection{Market profiles}
Autoregressive models are also commonly used to represent temporal dependence in both electricity and fuel prices \cite{weron2008forecasting}. Fuel price profiles, as well as power import and export prices, are generated using the same seasonal process as the aggregate heat demand profile without the subsequent spatial splitting step. The generated instances include up to six fuel types: natural gas, synthetic gas, biomass, biomethane, oil, and coal. For each fuel, monthly admissible price intervals are specified based on reported European district heating fuel cost ranges, reflecting realistic prices \cite{fuelprices}. Natural gas and synthetic gas are assigned wider intervals to capture stronger seasonal effects and higher short-term volatility, while biomass, biomethane, oil, and coal exhibit smoother variation. For electricity markets, prices can either be generated using the autoregressive model or sampled as a uniformly distributed ratio relative to a reference price level.

\subsection{Calibration against the Berlin reference network}

Although the procedure described in Section \ref{comp gen} - \ref{gen temp} is system-independent, it is calibrated against a specific real-world system to ensure that the generated instances remain realistic and transferable to proprietary datasets. In particular, the synthetic instances are calibrated to reproduce key operational and economic characteristics of real district heating systems, using a detailed reference model of the Berlin district heating network. Calibration targets structural indicators, including the number of production sites, the share of conversion technologies, and the typical ratio of installed heat capacity to peak demand. Temporal profiles, such as seasonal demand patterns and fuel price variability, are generated to match the statistical properties observed in the reference system. While the generated parameters are not exactly aligned with the reference values, their sampling ranges are derived from the Berlin instance and extended (e.g., by approximately $\pm 10 \%$) to introduce controlled variability around realistic baselines.

\section{Solution approaches} \label{sec:solution}
For the purpose of evaluating the generated dataset, we optimize a multi-objective unit commitment problem, which is modeled as a mixed-integer program (MIP). 

\subsection{MIP formulation}

Operation optimization of the multi-energy system is formulated as MIP, following the approach in \cite{Clarneretal2020}. The system consists of resources $R$, generating units $I$, storage units $K$, markets $M$, demand nodes $D$, and balance nodes $B$, represented by a network graph $G=(V,A)$ with $V = I \cup K \cup M \cup D \cup B$ and arcs $A$ defined by the network topology. All inputs are discretized over a set of time steps $T$. The MIP captures unit commitment, production, and network flow decision variables:
\begin{itemize}
    \item $z_{i,t} \in \{0,1\}$: unit status,
    \item $s_{i,t} \in \{0,1\}$: status changes,
    \item $x_{t,v}^{r} \ge 0$: incoming/outgoing flows of resource $r \in R$ at node $v \in V$,
    \item $h_{t,k}^{r}$: storage levels,
    \item $p_t^r, e_t^r \ge 0$: purchased and sold resources.
\end{itemize}

\resizebox{0.96\linewidth}{!}{
\begin{minipage}{1.03\linewidth}
\begin{align}
    &\text{min} && \left(\text{costs, emissions, $-$ CHP heat} \right) &&& \label{mip}\tag{$MIP$}\\
    & && \sum_{i \in I}  x_{t, i^{\text{out}}}^{r} + \sum_{k \in K} x_{t, k^{\text{out}}}^{r} + p_t^r  = d^r_t + \sum_{i \in I}  x_{t, i^{\text{in}}}^{r} + \sum_{k \in K} x_{t, k^{\text{in}}}^{r} + e_t^r &&& \forall r \in R, t \in T \label{eq:balance}\\
    & && x_{t, i^{\text{out}}}^{r_2} = \varphi_{i,t}^{r_1,r_2} \left(x_{t, i^{\text{in}}}^{r_1}\right) &&& \forall i \in I, t\in T, r_1, r_2 \in R \label{eq:conversion}\\
    & && s_{i,t} \leq z_{i,t}, s_{i,t} \leq 1 - z_{i,t}, s_{i,t} \geq z_{i,t} - z_{i,t-1} &&& \forall i \in I, t\in T \label{eq:activation}\\
    & && \sum_{\tau \in T^{\text{up}}_{i,t}} ( s_{i,t} - z_{i,\tau}) \leq 0, \sum_{\tau \in T^{\text{down}}_{i,t}} ( s_{i,t} + z_{i,\tau} -1) \leq 0 &&& \forall i \in I, t\in T \label{eq:minupdown}\\
    & && x_{t+1, i^{\text{out}}}^{r} - x_{t, i^{\text{out}}}^{r} \leq a_i^{\text{up}}, x_{t, i^{\text{out}}}^{r} - x_{t+1, i^{\text{out}}}^{r} \leq a_i^{\text{down}} &&&\forall i \in I, t\in T, r \in R \label{eq:ramping}\\
    & && h^r_{t+1,k} = a^{\text{loss}}_{t,k}h^r_{t,k} + a^{\text{load}}_{t,k}x_{t, k^{\text{in}}}^{r} - a^{\text{unload}}_{t,k}x_{t, k^{\text{out}}}^{r} &&&\forall k \in K, t\in T, r \in R \label{eq:storage}\\
    & && h,x,p,e \leq q_{\max}, \qquad h,x,p,e \geq q_{\min} \label{eq:capacities} \\
    \nonumber
\end{align}
\end{minipage}
}

The objective minimizes operational cost and CO$_2$ emissions while maximizing heat production from combined heat and power plants. Key constraints include demand and resource balance (\ref{eq:balance}), resource conversion via a piecewise-linear conversion map $\varphi$ (\ref{eq:conversion}), technical unit constraints (\ref{eq:activation}-\ref{eq:ramping}), storage dynamics (\ref{eq:storage}), and capacity limits (\ref{eq:capacities}). For a detailed description of the model, see \cite{riedmuller2025enhancing}.

\subsection{Algorithms}
\label{sec:solving}

Multi-objective optimization can be approached using a variety of methods, including weighted sum, lexicographic ordering, and the $\varepsilon$-constraint method \cite{Ehrgott2005}. Dynamic variants of the $\varepsilon$-constraint approach have previously been applied to district heating production portfolio planning \cite{Zittel2024}. To ensure reproducibility and consistent comparison across different instances, we adopt a lexicographic optimization strategy. Objectives are prioritized sequentially: operational cost is minimized first, followed by CO$_2$ emissions, and finally the negative heat production of combined heat and power units. In each stage, the solution of the preceding objective is treated as a constraint, preserving optimality with respect to higher-priority objectives. This approach ensures a unique and well-defined solution direction across instances.

\section{Computational study}
\label{sec:computational}

We evaluate the performance of the proposed instance generation and solution approach under a controlled computational setup. All experiments are performed using the Gurobi Optimizer~13.0.0 \cite{gurobi} (build v13.0.0rc1, Linux 64-bit, Debian GNU/Linux 12) on an Intel Xeon Gold 6338 CPU at 2.00\,GHz, with up to 32 threads utilized per run. Lexicographic optimization is applied with three objectives prioritized sequentially: operational cost, CO$_2$ emissions, and CHP heat production. Each objective is solved with a time limit of 10\,hours, resulting in a total maximum runtime of 30\,hours per instance. All solver parameters are kept at their default values, except for a MIP optimality gap of $10^{-6}$. To accelerate convergence, warm starts are used between consecutive objectives. The additional lexicographic constraint is enforced with a tolerance of $10^{-1}$ and scaled by $10^{-3}$ to improve numerical stability.

\subsection{Test instances}

The test instances considered in this section comprise the real-world reference model and the set of synthetically generated benchmarks.

\subsubsection{The real-world Berlin district heating model -- \texttt{B--VG2}}

Reference data of the real-world Berlin district heating network is provided by BEW Berliner Energie und Wärme GmbH to create instance \texttt{B--VG2}. The data involves the topology, parameters and time series data for the VG2 (Versorgungsgebiet 2) subgrid, which covers the East-Berlin section of the main interconnected network (Verbundnetz) in Berlin. For historical reasons, the district heating systems of East and West Berlin developed independently and were only later connected. Although today they are technically linked, the coupling remains limited, and the distinction between the two subgrids continues to play an operational role. 
The modeled subgrid supplies district heating to approximately 350,000 households. The time series dataset spans up to 25 years with a temporal resolution of 4 hours. The topology includes 5 substations with in total 5 combined heat and power plants (CHP), 2 of them combined cycle gas turbines, 5 heating plants, 5 power-to-heat plants (electric boilers), 6 heat pumps (heat sources: geothermal, data center, waste water, river). Thermal energy storage is not included in this instance, as storage facilities are currently only available in the West Berlin subgrid. The individual heat generation assets provide capacities of up to 924 MW of thermal output and up to 502 MW of electrical power. Most units are predominantly operated with gas \cite{BEW}.
The data further incorporates a single demand for the whole grid and market representations for electricity, CO$_2$ emissions, and multiple fuels, including natural gas, synthetic gas, oil, biomass, biomethane, and coal.

\subsubsection{Description of the synthetic instances -- \texttt{uc--[...]}}


The synthetic instances are based on the Berlin district heating network, serving as a conceptual reference for feasible topologies and component compositions. Building upon this reference, the instances are systematically generated with variations across all stochastic parameters, enabling the exploration of a diverse set of potential configurations:
\begin{itemize}
    \item length of the time horizon,
    \item number of demand nodes and production sites,
    \item number of converter units,
    \item inclusion or exclusion of storage units,
    \item number of fuel markets.
\end{itemize}

A total of 100 instances are generated. To ensure comparability, key global parameters are fixed across all instances, including the probability of transportation costs, the structure of network limiters, and the configuration of power markets, with exactly one import and one export market in each case. We generate instances with a temporal resolution of four hours and horizons of either 10 or 25 years, corresponding to 21\,900 and 54\,750 time steps, respectively. The instances are organized into groups sharing the same structural configuration, resulting in a largely consistent topology within each group. For each configuration, 10 stochastic realizations are generated to capture variability in converter characteristics, time series (e.g., demand and prices), and other operational inputs. Table~\ref{tab:instance_characteristics} summarizes the main characteristics of the benchmark instances. Each group (e.g., \texttt{uc\_000-uc\_009}) represents one configuration, while individual instances correspond to its stochastic realizations. The instances highlighted in gray serve as baseline references, both with and without storage. All other instances differ from one of these baselines by one or two parameters, which are highlighted in bold. The instance \texttt{B--VG2} represents the real-world Berlin reference model and serves as a baseline for comparison.

\begin{table}[t]
\centering
\caption{Features of the benchmark instances.}
\label{tab:instance_characteristics}
\begin{tabular}{lrrrrrr}
\hline
Instance
& $T$ 
& Demand 
& Prod.\ Sites 
& Converters
& Storage 
& Fuel Markets \\
\hline
B--VG2 & 54750 & 3 & 5 & 20 & 0 & 2 \\
\Xhline{0.1pt}
\rowcolor{gray!10}
uc\_000 to uc\_009 & 54750 & 3 & 5 & 20 & 0 & 2 \\
uc\_010 to uc\_019 & 54750 & \textbf{1} & 5 & 20 & 0 & 2 \\
uc\_020 to uc\_029 & 54750 & 3 & \textbf{3} & 20 & 0 & 2 \\
uc\_030 to uc\_039 & 54750 & 3 & 5 & \textbf{10} & 0 & 2 \\
uc\_040 to uc\_049 & \textbf{21900} & 3 & 5 & 20 & 0 & \textbf{4} \\
\Xhline{0.1pt}
\rowcolor{gray!10}
uc\_050 to uc\_059 & 54750 & 3 & 5 & 20 & 1 & 2 \\
uc\_060 to uc\_069 & 54750 & \textbf{5} & 5 & 20 & 1 & 2 \\
uc\_070 to uc\_079 & 54750 & 3 & 5 & 20 & 1 & \textbf{4} \\
uc\_080 to uc\_089 & 54750 & 3 & 5 & 20 & 1 & \textbf{6} \\
uc\_090 to uc\_099 & \textbf{21900} & 3 & 5 & 20 & 1 & 2 \\
\hline
\end{tabular}
\end{table}

\subsection{Data formats and availability}

The generated dataset of 100 instances is available on Zenodo (Project 19550701) \cite{Zenodo} in JSON format to facilitate use by energy modelers. In this format, each file encodes the topology of the multi-energy instance as a graph network, including its nodes, edges, associated parameters and time series. Additionally, MPS files are provided for direct use with optimization solvers; these files contain the full MIP formulation of each instance for each lexicographic 
optimization objective, along with the corresponding constraints.

\subsection{Computational results}
\begin{figure}[h!]
    \centering
    \includegraphics[scale=0.41]{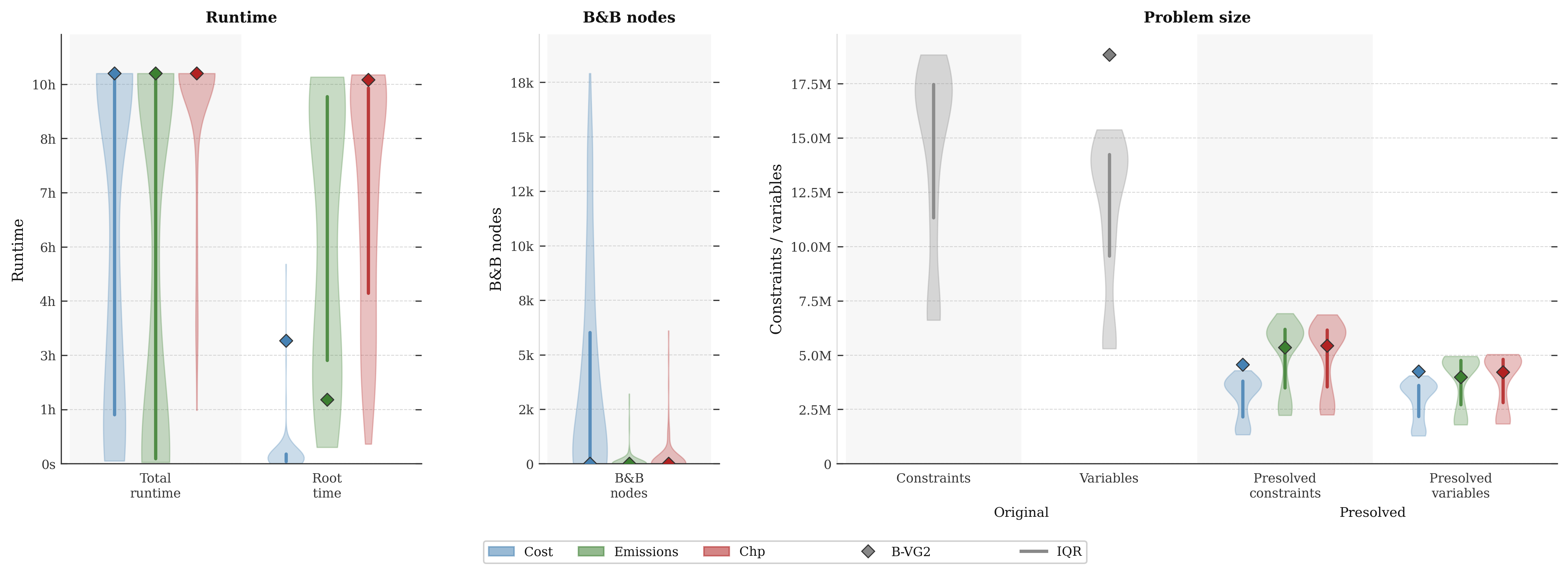}
    \caption{Distribution of characteristics over all instances for the three lexicographic objectives (1. cost, 2. emissions, 3. CHP).}
    \label{plot_violin}
\end{figure}

Fig.~\ref{plot_violin} illustrates the distribution of problem characteristics (e.g., number of variables, constraints, and branch-and-bound nodes) for the synthetic instances, with the real-world VG2 network shown as a diamond marker for reference. Presolve reductions are most effective for the first objective, while significantly fewer variables and constraints are eliminated for the second and third objectives. In lexicographic optimization, the use of a warm start when solving the MIP with an added constraint alters the problem structure; as a result, we can observe a reduction in the effectiveness of presolve. The ratio of discrete variables to the total number of variables (i.e., discrete plus continuous) remains stable across cases, but is lower for the second and third objectives. This indicates a higher share of continuous variables induced by the lexicographic constraint.

\begin{figure}[h!]

    \centering
    \includegraphics[scale=0.7]{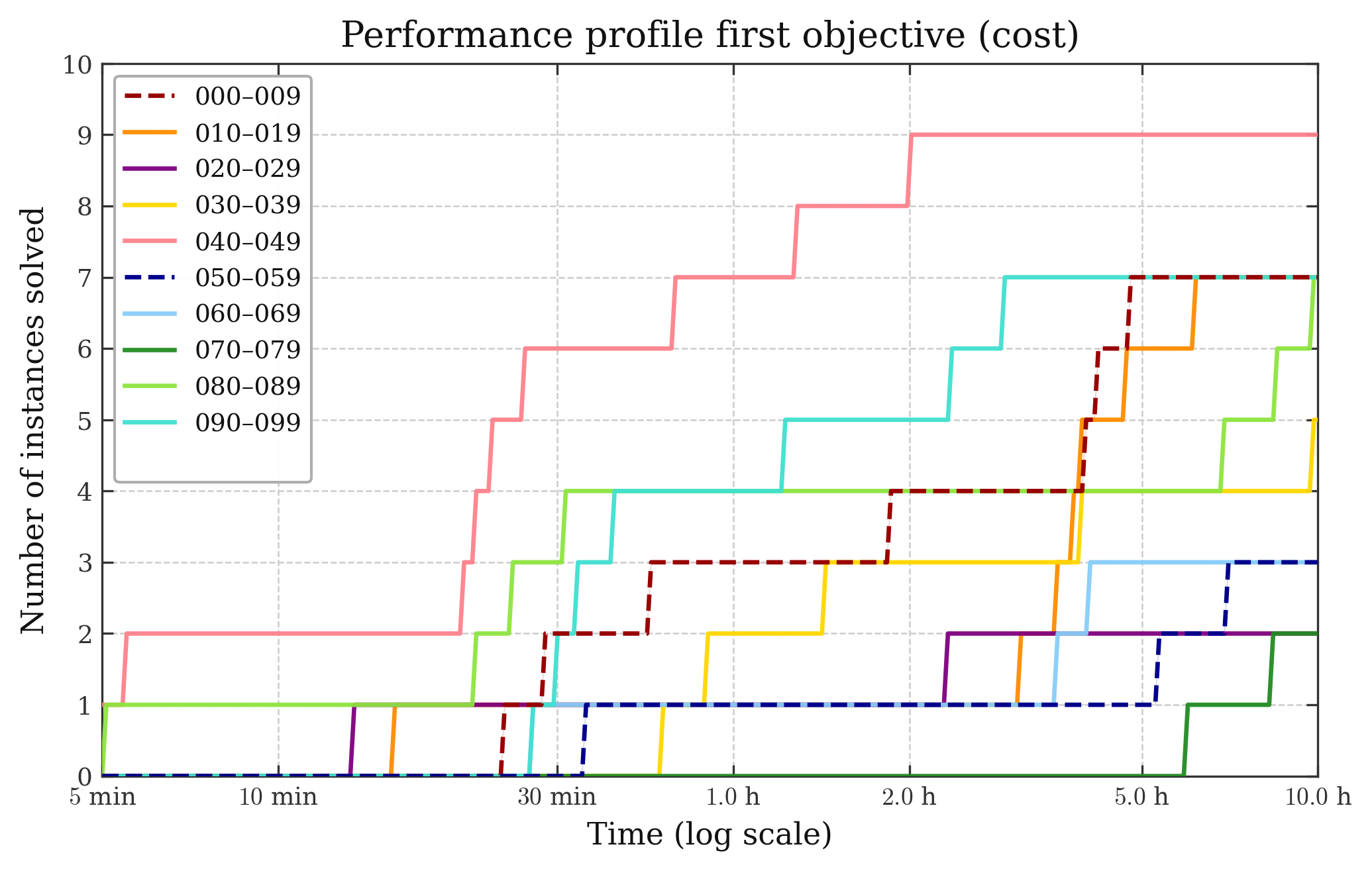}
    \caption{Runtime of first objective (cost) across each instance configuration.}
    \label{plot_cost}
\end{figure}

From a computational perspective, the first objective yields significantly larger branch-and-bound trees, highlighting its combinatorial complexity. In contrast, the second and third objectives require more root LP iterations and exhibit a wider RHS range, indicating increased numerical difficulty. This shift is also reflected in runtime: while a large fraction of instances for the first objective (cost minimization) is solved to optimality within the 10-hour limit, this fraction decreases markedly for subsequent objectives. Presolve time accounts for less than $1\%$ of total runtime across all instances. Scaling the additional constraint is essential, without it, performance deteriorates substantially. Overall, the results indicate a transition from predominantly combinatorial complexity in the first objective to increased numerical difficulty in the lexicographic extensions. Fig.~\ref{plot_cost} shows the number of solved instances over time on a logarithmic scale. The two baseline instance sets defined in Table~\ref{tab:instance_characteristics} are depicted as dashed lines and serve as reference datasets. The results indicate that
\begin{itemize}
    \item including storage increases computational effort: 30/50 instances without storage were solved within the time limit (mean over solved instances: 8,100 s), compared to 22/50 with storage (mean over solved instances: 13,745 s).
    
    \item instances with 10-year time series (\texttt{uc\_040-uc\_049}, \texttt{uc\_090-uc\_099}) are solved both faster and more consistently. This behavior can be attributed to the reduced number of variables and constraints. Overall, their average solving time decreases by 53.5\% relative to the baseline with 25-year time series.
    
    \item reducing the number of production sites to three (\texttt{uc\_020-uc\_029}) or limiting the number of converters to ten (\texttt{uc\_030-uc\_039}) increases computational difficulty.
    
    \item changes in the number of fuel markets have mixed effects: four markets (\texttt{uc\_070-uc\_079}) slightly reduce solving time compared to the two-market baseline, while six markets (\texttt{uc\_080}-\texttt{uc\_089}) increase it. Overall, these effects are small and influenced by time series and stochastic parameters.
    
    \item varying the number of demands shows minor and partially counterintuitive effects: five demands (\texttt{uc\_060-uc\_069}) marginally ease the problem, while a single demand (\texttt{uc\_010}-\texttt{uc\_019}) increases difficulty, yet their solving times remain close to baseline values.
\end{itemize}

 \subsubsection{Comparison with the Berlin model}
\label{sec:comparison}
We regard a synthetic instance as a realistic stand-in for the Berlin reference when
(i) its structural metrics (variable, constraint, and branch-and-bound-node counts) fall
within the reference's neighborhood, and (ii) its solver behavior is comparable. The no-storage
baseline group (\texttt{uc\_000-uc\_009}) satisfies both. It shares the VG2 topology,
differing mainly in temporal aspects. The number of variables and constraints is substantially higher due to a different representation of components (units connected in series, which
\texttt{B--VG2} models explicitly, are represented by a single converter). However, presolve eliminates structurally implied or non-binding constraints, such that the reduced model size of \texttt{B--VG2} falls inside the synthetic distribution (the VG2 diamond lies within the bulk of Fig.~\ref{plot_violin}). The
solver profile is likewise comparable: all three objectives remain unsolved within the 10-hour
limit, as for a subset of the synthetic instances. With both criteria met, these
observations support the conclusion that the synthetic instances provide a realistic
representation of the Berlin model.

\section{Conclusion and outlook}
\label{sec:conclusion}

We present a methodology for generating realistic synthetic benchmark instances of large-scale district heating systems, using a real-world model from Berlin as a reference. The framework is transferable to other cities and multi-energy networks. We provide three main contributions: a validated instance generation method, a computational study of factors affecting problem difficulty, and a dataset of 100 topology files with
corresponding optimization formulations, alongside the
open-source generator used to create them.
The realistic synthetic instances exhibit wide variability in computational difficulty, demonstrating their suitability for benchmarking. 
The results show that computational complexity depends on temporal, structural, and formulation factors; finer temporal resolution and storage increase model size and runtime, while structural simplifications may increase combinatorial complexity. The proposed framework is particularly valuable for assessing whether an existing model, which is already challenging to solve, exhibits increased computational difficulty under topological or stochastic modifications or whether such changes have a negligible effect on computational effort.

\section*{Nomenclature}
\label{sec:symbole_abk}

\subsection*{Abbreviations}
\begin{tabular}{@{}p{3cm}l}
    CHP & Combined heat and power plant \\
    MIP & Mixed-integer program \\
\end{tabular}

\subsection*{Latin Symbols}
\begin{tabular}{@{}p{3cm}l}
    $d, a, q, r_{\text{f}}, r_{\text{c}}$ & Pector for demand, coefficients, capacities, emission factor, emission certificate price \\
    $c,m, p, e$ & Variables for converters, markets, purchased resources and sold resources \\
    $G = (V,A)$ & Graph $G$ with set of vertices $V$ and set of arcs $A$ \\
    $R, I, K, M, D, B, C$ & Set of resources, generating units, storage units, markets, demands, balance nodes \\
    $c,m$ & Variables for converters, markets \\
    $T, S$ & Set of time steps and months \\
    $w, y$ & Weight of fuel assignment, source/ target series for conversion curve\\
    $x, z, s, h$ & Variables for resource flow, operation, status and storage level \\
\end{tabular}

\subsection*{Greek Symbols}
\begin{tabular}{@{}p{3cm}l}
    $\alpha, \beta, \rho$ & Ratios \\
    $\delta, \theta$ & Temperature level and fluctuation \\
    $\delta^+$ & Number of outgoing edges \\
    $\phi, \varepsilon$ & Parameter of autocorrelation and noise \\
    $\kappa, \lambda$ & Probabilities of capacity and cost nodes \\
    $\mu, \sigma$ & Temperature reference level and dispersion \\
    $\omega, \eta$ & Variable and noise of structural shares of the demand per district \\
    $\varphi, \vartheta$ & Resource conversion map \\
    $\xi$ & Amplitude factor \\
\end{tabular}

\subsection*{Superscripts and Subscripts}
\begin{tabular}{@{}p{3cm}l}
    $i/k$ & Generating/ storage unit \\
    $in/ out$ & Ingoing/ outgoing flow \\
    $r$ & Resource \\
    $t, s, k$ & Time step, monthly step, position of time step within a day \\
    $\mathcal{U}$ & Uniform distribution \\
\end{tabular}

\section*{Acknowledgments}
The work for this article has been conducted in the Research Campus MODAL funded by the German Federal Ministry of Research, Technology and Space (BMFTR) (fund numbers 05M14ZAM, 05M20ZBM, 05M2025).

\end{document}